\advance\year by -2000

\def\os{\obeyspaces}
\overfullrule=0pt\baselineskip=24pt
\def\C{{\rm C\hskip -6pt \vrule \hskip
6pt}}

\font\bigbold=cmbx17

\def\Prob{\hbox{Prob}}

\def\R{\os{\bf{R}\os}}

\def\C{\os{\bf{C}\os}}
\def\D{\os{\bf{D}\os}}

\def\S{\os{\bf{S}\os}}

\def\G{\os{\bf{G}\os}}
\def\H{\os{\bf{H}\os}}
\def\T{\os{\bf{T}\os}}

\def\NN{\os{\bf{NN}\os}}
\def\NF{\os{\bf{NF}\os}}
\def\FN{\os{\bf{FN}\os}}
\def\FF{\os{\bf{FF}\os}}
\def\F{\os{\bf{F}\os}}
\def\N{\os{\bf{N}\os}}
\def\HH{\os{\hbox{$\H\otimes\H$}}\os}
\def\HT{\os{\hbox{$\H\otimes\T$}}\os}
\def\TH{\os{\hbox{$\T\otimes\H$}}\os}
\def\TT{\os{\hbox{$\T\otimes\T$}}\os}

\def\R{{\bf R }}

\font\eightbf=cmr8

\def\p{\os{\bf{p}\os}}
\def\q{\os{\bf{q}\os}}

\line{\hfill
\the\month/\the\day/\the\year}

\centerline{\bigbold Quantum Game Theory}
\centerline{\bf by}
\centerline{\bf Steven E. Landsburg}
\centerline{\bf University of Rochester}

\centerline{to appear in}
\centerline{The Wiley Encyclopedia of Operations Research and Management 
Science}
\bigskip

{\it Quantum game theory} is the study of strategic 
behavior by agents with access to quantum technology.  
Broadly speaking, this technology can be employed in 
either of two ways:  As part of a randomization device 
or as part of a communications protocol.  

When it is used for randomization, quantum technology 
allows players to coordinate their strategies in certain 
ways.  The equilibria that result are all correlated 
equilibria in the sense of Aumann [A], but they form a 
particularly interesting subclass of correlated 
equilibria, namely those that are both achievable and 
deviation-proof when players have access to certain 
naturally defined technologies.  Not all correlated 
equilibria can be implemented via quantum strategies, 
and 
of those that can, not all are quantum-deviation-proof.

When players have access to private information, the 
theories of correlated and quantum equilibrium diverge 
still further, with the classical equivalence between 
mixed and behavioral strategies breaking down, and the
appearance of new equilibria that have no classical
counterparts.

The second game theoretic application of quantum 
technology, other than randomization, is to 
communication.  This leads to a new set of equilibria 
that seem to have no natural interpretation in terms of 
correlated equilibria or any other classical concepts. 

In Section I below, we will review the elements of game 
theory, with special emphasis on those aspects that are 
generalized or modified by quantum phenomena.  In 
Section II, we survey games with quantum randomization 
and in Section III, we survey games with quantum 
communication.

\centerline{\bf I.  Game Theory}

\noindent{\bf IA. Games.}

A {\it two-person game\/} consists of two sets $S_1$ and 
$S_2$ 
and two maps $$P_1:S_1\times S_2\rightarrow \R
\qquad\qquad
P_2:S_1\times S_2\rightarrow \R
$$
where \R denotes the real numbers.  The sets $S_i$ are 
called {\it strategy sets\/} and the functions $P_i$ are 
called {\it payoff functions\/}.  Games are intended to 
model strategic interactions between agents who are 
usually called {\it players\/} (though the players are 
not part of this formal definition).  The value 
$P_i(x,y)$ 
is called the {\it payoff\/} to Player $i$ when Player 1 
chooses strategy $x$ and Player 2 chooses strategy $y$.

For simplicity, we will usually assume the sets $S_i$ 
are
finite.

A {\it solution concept\/} is a function that associates 
to each game a subset of $S_1\times S_2$; the 
idea is to pick out those pairs of strategies that we 
believe players might select in real world situations 
modeled by the game.  The appropriate solution concept 
depends on the intended application.  The most studied 
solution concept is {\it Nash equilibrium}.  A pair 
$(x,y)$ is called a Nash equilibrium if $x$ maximizes 
the function $P_1(-,y)$ and $y$ maximizes the function 
$P_2(x,-)$.

\noindent{\bf IB.  Mixed Strategies}

To provide an accurate model of real world strategic 
situations, we must allow for the possibility that 
players might bend the rules.  For example, instead of 
choosing a single strategy, one or both players might 
randomize.  Starting with a game \G, we model this 
possibility by constructing the {\it associated mixed 
strategy game\/} $\G^{\hbox{\bf mixed}}$ in which the 
strategy space $S_i$ is replaced with the set 
$\Omega(S_i)$ of probability distributions on $S_i$, and 
the payoff function $P_i$ is replaced 
with the function  
$$P_i^{\hbox{mixed}}:(\mu,\nu)\mapsto\int 
P_i(x,y)d\mu(x)d\nu(y)$$

Although $\G^{\hbox{\bf mixed}}$ is not the same game as 
\G, it is traditional to refer to a Nash equilibrium in 
$\G^{\hbox{\bf mixed}}$ as a {\it mixed strategy 
equilibrium\/} in the game $\G$.   

An alternative but equivalent model would allow player 
$i$ to choose not a probability distribution on $S_i$ 
but an $S_i$-valued random variable from some allowable 
set.  This is the approach we will generalize in what 
follows.
  
\noindent{\bf IC.  Correlated Strategies}

In the play of $\G^{\hbox{\bf mixed}}$, we can imagine 
that Player $i$ first selects a probability 
distribution, then selects a random variable (with 
values in $S_i$) that realizes that distribution, 
then observes a realization of that random variable, and 
then plays accordingly.  Implicit in this description is 
that the random variables available to Player 1 are 
statistically independent of those available to Player 
2.

In the real world, however, this isn't always true.  In 
the extreme case, both agents might be able to observe 
the {\it same\/} random variable---such as the price of 
wheat as reported in the {\it Wall Street Journal\/}.  
We can model this extreme case by replacing $\G$ with a 
set of games $\{\G_\alpha\}$, one for each value 
$\alpha$ of the jointly observed random variable.  We 
then analyze each game $\G_\alpha$ separately.

But in less extreme cases, we need a new concept.  An {\it environment\/}
for {\bf G} is a pair $({\cal X}_1,{\cal X}_2)$, where ${\cal X}_i$ is 
a set of $S_i$-valued random variables.  (We do {\it not\/} assume that
random variables in ${\cal X}_1$ are necessarily independent of those in 
${\cal X}_2$.)  It is natural to assume that for any $X\in{\cal X}_i$ and
any map $\sigma:S_i\rightarrow S_i$, the random variable $\sigma\circ X$
is also in ${\cal X}_i$.  (This models the notion that players should be
able to map realizations to strategies any way they want to.)  

Given such
an environment, we define a new game ${\bf G}(E)={\bf G}({\cal X}_1,
{\cal X}_2)$ as follows:

Player $i$'s strategy set is ${\cal X}_i$.  The payoff functions are 
defined in the obvious way, namely

$$P_i(X,Y)=\int_{S_1\times S_2}P_i(x,y)d\mu_{X,Y}(x,y)$$
where $\mu_{X,Y}$ is the joint probability distribution on $S_1\times S_2$
induced by $(X,Y)$.

{\bf Definition IC.1.} The pair of random variables 
$(X,Y)$
is called a {\it correlated equilibrium\/} in \G if it 
is
a Nash equilibrium in the game  $\G(\{X\},\{Y\})$.  Two 
correlated equilibria are {\it equivalent\/} if they 
induce the
same probability distribution on $S_1\times S_2$.  We 
will
frequently abuse language by treating equivalent 
correlated
equilibria as if they were identical.

It is easy to prove the following:

{\bf Proposition IC.2.} Let $E$ be an environment and 
suppose
that $(X,Y)$ is a 
Nash equilibrium in the game
$\G(E)$.  Then $(X,Y)$ is
a correlated equilibrium in \G.

However. the converse to IC.2 does not hold:

{\bf Example IC.3.}  Let $S_1=S_2=\{\C,\D\}$.  Let
$X$, $Y$, and $W$ be random variables such that
$$\Prob(X=W=\C)=\Prob(X=W=\D)=1/8\qquad 
\Prob(X\neq W=\C)=\Prob(X\neq W=\D)=3/8$$
$$\Prob(Y=W=\C)=\Prob(Y=W=\D)=1/12\qquad \Prob(Y\neq 
W=\C)=\Prob(Y\neq W=\D)=
5/12$$
Let $E=(\{X,Y\},\{W\})$.
Let \G be the game with the following payoffs:

$$\matrix{
{\bf \hbox to 1.3in{}
\hbox{\bf Player Two}}\cr
{\vbox to 1in{\vfil \hbox{\bf Player
One}\vfil}
\hskip .2in\vbox {\offinterlineskip
\eightbf \halign{
\strut\hfil#&\quad\vrule#&\quad\hfil#
 \hfil&\quad\vrule#&
\quad\hfil#\hfil&\quad\strut\vrule#\cr
 &&{\bf C}&&{\bf D}&\cr
 \omit&&\omit&&\omit&\cr
 \noalign{\hrule}
\omit&&\omit&&\omit&\cr
 {\bf C}&&$(0,0)$&&$(2,1)$&\cr
 \omit&&\omit&&\omit&\cr
 \noalign{\hrule}
 \omit&&\omit&&\omit&\cr
 {\bf D}&&$(1,2)$&&$(0,0)$&\cr
 \omit&&\omit&&\omit&\cr
 \noalign{\hrule}
 }}}\cr}$$
\vskip 4pt\noindent

It is easy to check that both $(X,W)$ and $(Y,W)$ yield
correlated equilibria in \G.    
But $(X,W)$ is not an equilibrium in the game $\G(E)$, 
though $(Y,W)$ is.

\noindent{\bf 1D.  Games with Private Information.}

In real world strategic interactions, either player 
might know something the other doesn't.   
We model this situation as a {\it game of private 
information\/}, consisting of two strategy spaces $S_i$, 
two sets (called {\it information sets\/} ${\cal A}_i$, 
a probability distribution on 
${\cal A}_1\times{\cal A}_2$, and two payoff functions
$$P_i:{\cal A}_1\times{\cal A}_2\times S_1\times 
S_2\rightarrow \R$$

Given a game of private information, the {\it associated 
game\/} $\G^\#$ has strategy sets ${S}_i^{\#}=Hom({\cal 
A}_i,{S}_i)$ and payoff functions
$$P_i^\#(F_1,F_2)=\int_{{\cal A}_1\times{\cal 
A}_2}P_i(A_1,A_2, F_1(A_1),F_2(A_2))$$

(Here $Hom(A,S)$ denotes the set of all functions from $A$ 
to $S$.)

If \G is a game of private information, a {\it Nash 
equilibrium\/} in \G is (by definition) a Nash 
equilibrium in the ordinary game $\G^{\#}$.

Now we want to enrich the model so players can 
randomize.  To this end, 
let $E$ be an environment in the sense of Section 1C; 
that is, $E=({\cal X}_1,{\cal X}_2)$ where ${\cal X}_i$ 
is a set of 
$S_i$-valued random variables.  We define the {\it 
associated environment\/} $E^{\#}=({\cal X}_1,{\cal 
X}_2)$ by setting 
${\cal X}_i^{\#}=Hom(A_i,{\cal X}_i)$ and identifying 
the latter set with a set of $S_i^{\#}$-valued random 
variables.  

Thus, if \G is a game of private information with 
environment $E$, we can first ``eliminate the private 
information'' by passing to the associated game 
$\G^{\#}$ and environment $E^{\#}$, and then ``eliminate 
the random variables'' by passing to the associated game 
$\G^{\#}(E^{\#})$ as in the discussion preceding 
Definition 1C.1.  Now we're studying an ordinary game, 
where we have the ordinary notion of Nash equilibrium.

Unfortunately, {\it this construction does not 
generalize to the quantum context\/}.  To get a 
construction that generalizes, we need to proceed in the 
opposite order, by first eliminating the random 
variables and then eliminating the private information:

{\bf Construction ID1.}  Given a game of private information
{\bf G} with an environment $E=({\cal X}_1,{\cal X}_2)$, define a new
game of private information ${\bf G}(E)$ as follows:

The information sets ${\cal A}_i$ and the probability distribution on
${\cal A}_1\times{\cal A}_2$ are as in {\bf G}.  The strategy sets are
the ${\cal X}_i$.  The payoff functions are

$$P_i(A_1,A_2,X,Y)
=\int_{S_1\times S_2}P_i(A_1,A_2,s_1,s_2)d\mu_{X,Y}(s_1,s_2)$$

Now applying the $\#$ construction to ${\bf G}(E)$ gives an ordinary 
game ${\bf G}(E)^\#$.

{\bf Theorem.}  $\G(E)^{\#}=\G^{\#}(E^{\#})$.  

{\bf Remark}.  In spirit, and in the language of [K], 
$\G^{\#}(E^{\#})$ is like the associated game with {\it 
mixed strategies\/}, where player $i$ chooses a 
probability distribution over maps $A_i\rightarrow S_i$, 
while $\G(E)^{\#}$ is like the associated game with {\it 
behavioral strategies\/}, where player $i$ chooses, for 
each element of $A_i$, a probability distribution over 
$S_i$.  In [K], these games are equivalent in an 
appropriate sense; here they are actually the same game.  
The difference is that in [K], players choose 
probability distributions whereas here they choose 
random variables.  But the two approaches are 
fundamentally equivalent.

\noindent{\bf IE.  Quantum Game Theory}

The notions of mixed strategy and correlated equilibria 
are meant to model the real-world behavior of strategic 
agents who have access to randomizing technologies (such 
as weighted coins).  {\it Quantum game theory\/} is the 
analogous attempt to model the behavior of strategic 
agents who have access to quantum technologies (such as 
entangled particles).  

Broadly speaking, there players might use these quantum 
technologies in either of two ways:  As randomizing 
devices or as communication devices.  We will consider 
each in turn.

\medskip

\centerline{\bf II.  Quantum Randomization}

{\bf IIA. Quantum Strategies.}  Just as the theory of 
mixed
strategies models the behavior of players with access to 
classical randomizing devices, the theory of quantum 
strategies
models the behavior of players with access to quantum 
randomizing
devices.  These devices provide players with access to 
families 
of observable quantities that cannot be modeled as 
classical 
random variables.

For example, let $X$, $Y$, $Z$ and $W$ be binary random 
variables.  
Classically we have the near triviality: 
   $$Prob(X\neq W)\le Prob(X\neq Y)+Prob(Y\neq 
    Z)+Prob(Z\neq W)$$
(Proof:  Imagine $X,Y,Z,W$ lined up in a row; in order 
for $X$ to differ from $W$, at least one of $X,Y,Z,W$ 
must differ from its neighbor.)  But if $X,Y,Z$ and $W$ 
are quantum mechanical measurements, this inequality 
need 
not hold.  (The most obvious paradoxes are avoided by 
the fact that neither $X$ and $Z$, nor $Y$ and $W$, are 
simultaneously observable.)

To model the behavior of agents who can make such 
 measurements, we can  mimic the definitions (I.C.1) and
 (I.C.2),  replacing the sets ${\cal X}_i$  
of random variables with sets ${\cal X}_i$  
  of quantum mechanical observables.    
 We require that any $X\in{\cal X}_1$ and any $Y\in{\cal 
X}_2$ 
be simultaneously observable.  We call such a pair 
$E=({\cal X}_1,
{\cal X}_2)$ a {\it quantum environment\/}.  (Quantum 
environments will be defined more precisely in Section 
II.B below.)  Given such a quantum environment and given 
a game \G, we  
construct a new game $\G(E)$ just as in the discussion 
preceding
 Definition IC.1.

If $(X,Y)$ is a Nash equilibrium in $\G(E)$, we 
sometimes speak loosely 
enough to say that $(X,Y)$ is a {\it quantum 
equilibrium\/} in \G, though
the property of being a quantum equilibrium depends not 
just on \G but 
on the environment $E$.  Two quantum equilibria are 
called {\it equivalent\/}
if they induce the same probability distribution on 
$S_1\times S_2$, and,
as with correlated equilibria, we will sometimes speak 
of equivalent 
quantum equilibria as if they were identical.

If $(X,Y)$ is a quantum equilibrium then (by the 
definition of
quantum environment), $X$ and $Y$ are simultaneously 
observable
and hence can be treated as classical random variables.  
With 
this identification, it is easy to show that any quantum 
equilibrium 
is a correlated equilibrium.  (This generalizes 
Proposition IC.2.)  
However, just as in Example IC.3, the converse need not 
hold.  A pair 
$(X,Y)$ that is an equilibrium in one quantum 
environment need not be
an equilibrium in another.

\noindent{\bf IIB.  The Quantum Environment}

Consider a game in which the strategy sets are 
$S_1=S_2=\{\C,\D\}$.  Players can implement (ordinary 
classical) mixed strategies by flipping (weighted) 
pennies, mapping the outcome ``heads'' to the strategy 
\C and the outcome ``tails'' to the strategy \D.

While a classical penny occupies either the state \H 
(heads up) or \T (tails up), a quantum penny can occupy 
any state of the form $\psi=\alpha \H+\beta \T$, where 
$\alpha$ and $\beta$ are complex scalars, not both zero.  
A heads/tails measurement of such a penny yields the 
outcome either \H or \T with probabilities proportional 
to $|\alpha|^2$ and $|\beta|^2$.  Physical actions such 
as rotating the penny induce unitary transformations of 
the state space, so that the state $\psi$ is replaced by 
$U\psi$ where $U$ is some unitary operator on the 
complex vector space $\C^2$.
 
(Of course literal macroscopic pennies do not behave 
this way, but spin-1/2 particles such as electrons do, 
with ``heads'' and ``tails'' replaced by ``spin up'' and 
``spin down''.)  

A single quantum penny is no more or less useful than a 
classical randomizing device.  If you want to play heads 
with probability $p$, you can first apply a unitary 
transformation that converts the state to some 
$\gamma\H+\delta\T$ with 
$|\gamma|^2/(|\gamma|^2+|\delta^2|)=p$, then measure the 
heads/tails state of the penny and play accordingly.

However, two players equipped with quantum pennies have 
something more than a classical randomizing device.  A 
pair of quantum pennies occupies a state of the form
$$\alpha\HH+\beta\HT+\gamma\TH+\delta\TT$$
where $\alpha, \beta,\gamma,\delta$ are complex scalars, 
not all zero.  A physical manipulation of the first 
penny transforms the first factor unitarily and a 
physical manipulation of the second penny transforms the 
second factor unitarily.  Subsequent measurements yield 
the outcomes (heads,heads), (heads,tails) and so forth 
with probabilities proportional to $|\alpha|^2$, 
$|\beta|^2$ and so forth.  This allows the players to 
achieve joint probability distributions that cannot be 
achieved via the observations of independent random 
variables.

{\bf Example II.B.1.}  Suppose that two pennies begin in 
the {\it maximally entangled state\/} $\HH+\TT$.  
Players One and Two apply the transformations $U$ and 
$V$ to the first and second pennies where
$$U=\pmatrix{\cos(\theta)&\sin(\theta)\cr-
\sin(\theta)&\cos(\theta)\cr}\qquad
V=\pmatrix{\cos(\phi)&\sin(\phi)\cr-
\sin(\phi)&\cos(\phi)\cr}$$ 
This converts the state from $\HH+\TT$ to 
$$U\H\otimes V\H+U\T\otimes V\T=\cos(\theta-\phi)\HH
+\sin(\theta-\phi)\HT-\sin(\theta-\phi)\TH+\cos(\theta-
\phi)\TT$$

If players map the outcomes \H and \T to the strategies 
\C and \D, then the resulting probability distribution 
over strategy pairs is 
$$\Prob(\C,\C)=\Prob(\D,\D)=\cos^2(\theta-\phi)/2\qquad
\Prob(\C,\D)=\Prob(\D,\C)=\sin^2(\theta-\phi)/2$$

We can now make precise the notion of {\it quantum 
environment}; a quantum enviroment is a triple 
$(\xi,{\cal X}_1,{\cal X}_2)$ where  
\itemitem{a)}$\xi$ is a non-zero vector in a complex 
vector space $\C^{n_1}\otimes \C^{n_2}$ (with $n_1$ and 
$n_2$ assumed finite here, though this could all be 
generalized) 
\itemitem{b)}${\cal X}_i$ is a set of unitary operators 
on $\C^{n_i}$.

The unitary operators fill the same role as the sets of 
random variables in Section I; they are the things that 
players can observe, and on whose realizations they can 
condition their strategies.  

\noindent{\bf IIC. Quantum Equilibrium}

Consider again the game \G from Example IC.3:

$$\matrix{
{\bf \hbox to 1.3in{}
\hbox{\bf Player Two}}\cr
{\vbox to 1in{\vfil \hbox{\bf Player
One}\vfil}
\hskip .2in\vbox {\offinterlineskip
\eightbf \halign{
\strut\hfil#&\quad\vrule#&\quad\hfil#
 \hfil&\quad\vrule#&
\quad\hfil#\hfil&\quad\strut\vrule#\cr
 &&{\bf C}&&{\bf D}&\cr
 \omit&&\omit&&\omit&\cr
 \noalign{\hrule}
\omit&&\omit&&\omit&\cr
 {\bf C}&&$(0,0)$&&$(2,1)$&\cr
 \omit&&\omit&&\omit&\cr
 \noalign{\hrule}
 \omit&&\omit&&\omit&\cr
 {\bf D}&&$(1,2)$&&$(0,0)$&\cr
 \omit&&\omit&&\omit&\cr
 \noalign{\hrule}
 }}}\cr}$$
\vskip 4pt\noindent

Let $E$ be the following quantum environment:  
$\xi=\HH+\TT$ where $\{\H,\T\}$ is a basis for the 
complex vector space $\C^2$;
${\cal X}_1={\cal X}_2$ is the set of all operators that 
take the form 
$$M(\theta)=\pmatrix{\cos(\theta)&\sin(\theta)\cr -
\sin(\theta)&\cos(\theta)\cr}$$
when expressed in terms of the basis $\{\H,\T\}$.

Physically, this means that each player has access to 
one member of a pair of maximally entangled pennies, and 
can rotate that penny through any angle before measuring 
it's heads/tails orientation.  

Taking Player Two's angle $\phi$ as given, Player One 
clearly optimizes by choosing $\theta$ so that 
$\sin(\theta-\phi)=1$, and symmetrically for Player One.  
Thus in equilibium, the outcomes (2,1) and (1,2) are 
each realized with probability 1/2.  

As we've noted earlier, this is of course a correlated 
equilibrium, but it is more than that.  For example, the 
correlated equilibria $(X,W)$ and $(Y,W)$ of Example 
IC.3 are not sustainable as quantum equilibria in this 
environment.

\noindent{\bf IID. Quantum Games of Private Information}

Let \G be a game of private information as defined in 
Section ID, and let $E$ be a quantum environment. 

We would like to model the play of \G in the environment 
$E$ as the play of an ordinary game.  A natural attempt 
is to replace \G with the game $\G^{\#}$ defined in 
Section ID; recall that a strategy in $\G^{\#}$ is a map 
from the information set ${\cal A}_i$ to the strategy 
set $S_i$.  However, in the quantum case there is no 
natural way to define an environment $E^{\#}$ for this 
game.  More precisely, it is shown in [DL] that the 
natural definition of $E^{\#}$ makes sense when and only 
when the measurements in ${\cal X}_1\cup{\cal X}_2$ have 
a single joint probability distribution, so that they 
can be thought of as classical random variables.  In 
other words, the existence of an environment $E^{\#}$ is 
equivalent to the absence of any specifically quantum 
phenomena.

Therefore we must generalize not the construction 
$\G^{\#}(E^{\#})$ from Section ID, but the classically 
equivalent construction $\G(E)^{\#}$.  In the language 
of [K] (and of the remark at the end of Section ID), 
quantum game theory allows players to choose behavioral 
strategies without equivalent mixed strategies.

{\bf Example IID.1.}   
Let \G be the following game of private 
information:

The information sets are ${\cal A}_1={\cal 
A}_2=\{\hbox{red,green}\}$.  The probability 
distribution on ${\cal A}_1\times{\cal A}_2$ is uniform. 
The payoff functions are specified as follows:  

$$\matrix{
{\bf \hbox to 1.3in{IF BOTH PLAYERS OBSERVE RED}}\cr
\hbox{\hfil{\bf \hbox{\bf 
\phantom{xxxxxxxxxxxxxxxx}Player Two}}\hfil}\cr
{\vbox to 1in{\vfil \hbox{\bf Player
One}\vfil}
\hskip .2in\vbox {\offinterlineskip
\eightbf \halign{
\strut\hfil#&\quad\vrule#&\quad\hfil#
 \hfil&\quad\vrule#&
\quad\hfil#\hfil&\quad\strut\vrule#\cr
 &&{\bf C}&&{\bf D}&\cr
 \omit&&\omit&&\omit&\cr
 \noalign{\hrule}
\omit&&\omit&&\omit&\cr
 {\bf C}&&$(1,1)$&&$(0,0)$&\cr
 \omit&&\omit&&\omit&\cr
 \noalign{\hrule}
 \omit&&\omit&&\omit&\cr
 {\bf D}&&$(0,0)$&&$(1,1)$&\cr
 \omit&&\omit&&\omit&\cr
 \noalign{\hrule}
 }}}\cr}\qquad
\matrix{
{\bf \hbox to 1.3in{IF EITHER PLAYER OBSERVES GREEN}}\cr
\hbox{\hfil{\bf \hbox{\bf 
\phantom{xxxxxxxxxxxxxxxx}Player Two}}\hfil}\cr
{\vbox to 1in{\vfil \hbox{\bf Player
One}\vfil}
\hskip .2in\vbox {\offinterlineskip
\eightbf \halign{
\strut\hfil#&\quad\vrule#&\quad\hfil#
 \hfil&\quad\vrule#&
\quad\hfil#\hfil&\quad\strut\vrule#\cr
 &&{\bf C}&&{\bf D}&\cr
 \omit&&\omit&&\omit&\cr
 \noalign{\hrule}
\omit&&\omit&&\omit&\cr
 {\bf C}&&$(0,0)$&&$(1,1)$&\cr
 \omit&&\omit&&\omit&\cr
 \noalign{\hrule}
 \omit&&\omit&&\omit&\cr
 {\bf D}&&$(1,1)$&&$(0,0)$&\cr
 \omit&&\omit&&\omit&\cr
 \noalign{\hrule}
 }}}\cr}$$
\vskip 4pt\noindent

The environment $E$ is as in IIC.  

In the game $\G(E)^{\#}$, players choose maps ${\cal 
A}_i\rightarrow {\cal X}_i$.  

To find equilibria in this game, suppose that Players 
One has chosen to map ``red'' and ``green'' to 
$M(\theta_{red})$ and $M(\theta_{green})$, while Player 
Two has chosen $M(\phi_{red})$ and M($\phi_{green})$ 
(where the M matrices are as defined in Section IC.)  
Then we can compute Player One's expected payoffs as functions on 
${\cal A}_1\times {\cal A}_2$:

$$\eqalignno{
EP_1(\hbox{red},\hbox{red})&=\cos^2(\theta_{red}-
\phi_{red})&(IID.1)\cr
EP_1(\hbox{red},\hbox{green})&=\sin^2(\theta_{red}-
\phi_{green})&(IID.2)\cr
EP_1(\hbox{green},\hbox{red})&=\sin^2(\theta_{green}-
\phi_{red})&(IID.3)\cr
EP_1(\hbox{green},\hbox{green})&=\sin^2(\theta_{green}-
\phi_{green})&(IID.4)\cr}$$

Because the probability distribution on ${\cal 
A}_1\times{\cal A}_2$ is uniform, Player One seeks to 
maximize the sum of these four expressions, taking 
$\phi_{red}$ and $\phi_{green}$ as given.  At the same 
time (due to the symmetry of the  problem) Player Two 
seeks to maximize the identical sum, taking 
$\theta_{red}$ and $\theta_{green}$ as given.

Given this, we can compute that there are two types of 
equilibria:

\itemitem{a)}Equilibria in which exactly three of the 
four expressions (IID.1)-(IID.4) are equal to 1 and the 
fourth is equal to 0.  All four possibilities occur with 
$\phi_r,\phi_g,\theta_r,\theta_g\in\{0,{\pi\over2},\pi\}
$.
In these equilibria each player receives a payoff of 
$3/4=.75$.

\itemitem{b)}Equilibria equivalent to 
$\{\phi_{red}=0,\phi_{green}=3\pi/4,\theta_{red}=\pi/8,
\theta_{green}=3\pi/8\}$.  In these equilibria each 
player receives a payoff of 
$(1/2+\sqrt{2}/4)\approx.85.$ 

The best equilibrium that can be reached in any 
classical environment is equivalent to an equilibrium of 
type a).  This is so even when the classical environment 
includes correlated random variables.  

For ordinary games, we observed that every quantum 
equilibrium is also a correlated equilibrium.  For games 
of private information, this example demonstrates that 
no analogous statement is true.

\medskip

\centerline{\bf III. Quantum Communication}

\noindent {\bf IIIA.  Communication}

In any real world implementation of a game, players must 
somehow communicate their strategies before they can 
receive payoffs.  This follows from the fact that the 
payoff functions take both players' strategies as 
arguments; therefore information about both strategies 
must somehow be present at the same place and time.  

Ordinarily, the communication process is left unmodeled.  
But in this section, we need an explicit model so that we can 
explore the effects of allowing quantum communication 
technology.  To that end we postulate a referee who 
communicates with the players by handing them markers 
(say pennies) which the players can transform from one 
state to another (say by flipping them over or leaving 
them unflipped) to indicated their strategy choices; the 
markers are eventually returned to the referee who 
examines them and computes payoffs.

Of course not all real world games have literal 
referees; sometimes the players are firms, their 
strategies are prices, and the prices are communicated 
not to a referee but to a marketplace, where payoffs are 
``computed'' via market processes.  
In the models to follow, the referee can be understood 
as a metaphor for such processes.

\noindent{\bf IIIB.  A Quantum Move}  

Consider the game with the following payoff matrix (for 
now, view the labels \NN, \NF, etc., as arbitrary labels 
for strategies):

$$\matrix{
{\bf \hbox to 1.3in{}
\hbox{\bf Player Two}}\cr
{\vbox to 1in{\vfil \hbox{\bf Player
One}\vfil}
\hskip .2in\vbox {\offinterlineskip
\eightbf \halign{
\strut\hfil#&\quad\vrule#&\quad\hfil#
 \hfil&\quad\vrule#&
\quad\hfil#\hfil&\quad\strut\vrule#\cr
 &&{\bf N}&&{\bf F}&\cr
 \omit&&\omit&&\omit&\cr
 \noalign{\hrule}
 {\bf NN}&&$(1,0)$&&$(0,1)$&\cr
 \noalign{\hrule}
 {\bf NF}&&$(0,1)$&&$(1,0)$&\cr
 \noalign{\hrule}
{\bf FN}&&$(0,1)$&&$(1,0)$&\cr
 \noalign{\hrule}
{\bf FF}&&$(1,0)$&&$(0,1)$&\cr
 \noalign{\hrule}
 }}}\cr}$$
\vskip 4pt\noindent

As noted above, classical game theory does not ask how 
players communicate with the referee.  If we want to 
embellish our model with an explicit communication 
protocol, there are several essentially equivalent ways 
to do it.  

{\bf IIIB.1.  The Simplest Protocol:}  The referee 
passes 
two pennies to Player One, who flips both to indicate a 
play of \FF, flips only the first to indicate a play of 
\FN, and so forth, and one penny to Player Two, who 
flips (\F) or does not (\N).  The pennies are returned 
to the referee, who examines their states and makes 
payoffs accordingly.  

{\bf IIIB.2.  An Alternative Protocol:}  A single penny 
in state \H is passed to Player One, who either flips or 
doesn't; the penny is then passed to Player Two, who 
either flips or doesn't; the penny is then returned to 
the Player One, who either flips or doesn't; the penny 
is then returned to the referee who makes payoffs of 
(1,0) if the final state is \H or (0,1) if the final 
state is \T.  (The players are blindfolded and cannot 
observe each others' plays.)  

We can set things up so that flipping and not-flipping 
correspond to the applications of the unitary matrices
$$\F=\pmatrix{0&1\cr -i&0\cr}\qquad \N=\pmatrix{1&0\cr 
0&1\cr}$$

Now suppose that Player One manages to cheat by 
employing arbitrary unitary operations. Under the 
simplest protocol IIIB.1, this is equivalent to playing 
a mixed strategy and gives Player One no advantage.  But 
under the alternative protocol, Player One can guarantee 
himself a win by choosing the unitary matrix
$$U={1\over2}\pmatrix{1+i&\sqrt{2}\cr -\sqrt{2}&
1-i\cr}$$ on  his first turn and $U^{-1}$ on his second.  
This guarantees him a win because $U^{-1}\circ F\circ 
U$ and $U^{-1}\circ N\circ U$ are both diagonal 
matrices, so that both fix the state \H.  (Recall
that states are unchanged by scalar multiplication.) In 
other 
words, the final state is now \H regardless of Player 
Two's strategy.  

Note that it might be quite impossible to prohibit 
Player One from employing the matrix $U$, for exactly 
the same reason that it is usually quite impossible to 
prohibit players for adopting mixed strategies:  When 
the referee makes his final measurement, the only 
information revealed is the final state of the penny, 
not the process by which it achieved that state.

This example, due to David Meyer ([M]), illustrates two 
points:  First, quantum communication can matter.  
Second, whether or not quantum communication matters 
depends not just on the game \G; it depends on the 
specified communications protocol.  

\noindent{\bf IIIC. The Eisert-Wilkens-Lewenstein 
Protocol}

The Eisert-Wilkens-Lewenstein protocol ([EW],[EWL]) 
captures the 
potential effects of quantum communication in a quite 
general context and is therefore the most studied 
protocol in games of quantum communication.  

We start with a game \G in which the strategy sets are 
$S_1=S_2=\{\C,\D\}$.  (Everything can be generalized to 
larger strategy sets, but we will stick to this simplest 
case.)  The referee prepares a pair of pennies in the 
state \HH+\TT, and passes one penny to each player, with 
Player One receiving the ``left-hand'' penny.  

Players are instructed to play operate with one of the 
unitary matrices \F and \N of Section IIIB, with \F indicating
a desire to play \D and \N a desire to play \C.  Players
actually operate with the unitary matrices of their choices
and then return the pennies to the referee, who makes a 
measurement that distinguishes among the four states that
could result if the players followed instructions.

We can of course model this situation by saying that the 
original game \G has been replaced by the {\it 
associated quantum game\/} $\G^Q$ (not to be confused 
with the associated quantum games that arise in the very 
different context of Section II), with strategy sets 
equal to the unitary group $U_2$ (or, equivalently---
because states are unchanged by scalar multiplication---
the special unitary group $SU_2$) and payoff functions
have the obvious definition.    It is observed in 
[L] that $\G^Q$ is equivalent to the following game:
\itemitem{$\bullet$}The strategy sets $S_i^Q$ are both 
equal to the group of unit quaternions (which is 
isomorphic to the special unitary group $SU_2$)
\itemitem{$\bullet$}The payoff functions are defined as 
follows:
$$P_i^Q(\p,\q)=A^2 P_i(\C,\C)+B^2 P_i(\C,\D)+
C^2 P_i(\D,\C)+D^2 P_i(\D,\D)$$
where $P_i$ is the payoff function in \G and where 
$$\p\q=A+Bi+Cj+Dk$$

The group structure on the strategy sets guarantees that 
(except in the uninteresting case where both payoff 
functions are maximized at the same arguments) there can 
be no pure-strategy equilibria in this game, because 
Player One, taking Player Two's strategy \q as given, 
can always choose \p to maximize his own payoff, in 
which case Player Two cannot be maximizing.  So the game 
$\G^Q$ is quite uninteresting unless we allow mixed 
strategies.

A mixed strategy in $\G^Q$ is a probability measure on 
the strategy space $S_1^Q=SU_2=\S^3$ where $\S^3$ is the 
three-sphere.  Thus the strategy spaces in $(\G^Q)^{\bf 
mixed}$ are quite large.  However, it is shown in [L] 
that in equilibrium, both players can be assumed to 
adopt strategies supported on at most four points, and 
that each set of four points must lie in one of a small 
number of highly restrictive geometric configurations.  
This considerably eases the problem of searching for 
equilibria.  

{\bf Example IIIC.1.  The Prisoner's Dilemma.}  Consider 
the ``Prisoner's Dilemma'' game

$$\matrix{
{\bf \hbox to 1.3in{}
\hbox{\bf Player Two}}\cr
{\vbox to 1in{\vfil \hbox{\bf Player
One}\vfil}
\hskip .2in\vbox {\offinterlineskip
\eightbf \halign{
\strut\hfil#&\quad\vrule#&\quad\hfil#
 \hfil&\quad\vrule#&
\quad\hfil#\hfil&\quad\strut\vrule#\cr
 &&{\bf C}&&{\bf D}&\cr
 \omit&&\omit&&\omit&\cr
 \noalign{\hrule}
\omit&&\omit&&\omit&\cr
 {\bf C}&&$(3,3)$&&$(0,5)$&\cr
 \omit&&\omit&&\omit&\cr
 \noalign{\hrule}
 \omit&&\omit&&\omit&\cr
 {\bf D}&&$(5,0)$&&$(1,1)$&\cr
 \omit&&\omit&&\omit&\cr
 \noalign{\hrule}
 }}}\cr}$$
\vskip 4pt\noindent

There is only one Nash equilibrium, only one (classical) 
mixed strategy equilibrium, and only one (classical) 
correlated equilibrium, namely $(\D,\D)$ in every case.   
But it is an easy exercise to check that there are 
multiple mixed-strategy equilibria in the  
Eisert-Wilkens-Lewenstein game $\G^Q$.

First example: Each player chooses any
four orthogonal quaternions and plays each of the four
with equal probability.  This induces the probability 
distribution in which each of the four payoffs is equiprobable
and each player earns an expected payoff of 9/4.

Second example:
Player One plays the quaternions $1$ and $i$ with equal 
probability and Player Two plays the quaternions $j$ and $k$
with equal probability.  This induces the probability 
distribution where  the payoffs 
$(0,5)$ and $(5,0)$ each occur with probability $1/2$ 
so that each player earns an expected payoff of 5/2.

The techniques of [L] reveal that, up to a 
suitable notion of equivalence, these are the 
only mixed strategy 
equilibria in $\G^Q$.

\vfill\eject

\centerline{\bf References}

\item{[A]}Aumann, ``Subjectivity and Correlation in
Randomized Strategies'', {\it J. Math Econ\/} 1 (1974).

\item{[CHTW]}R.~Cleve, P.~Hoyer, B.~Toner, and
J.~Watrous,
``Consequences and Limits of Nonlocal Strategies'', {\it
Proc. of the 19th Annual Conference on Computational
Complexity\/} (2004), 236-249

\item{[DL]}G.~Dahl and S.~Landsburg, ``Quantum 
Strategies'', preprint 
available at http://www.landsburg.com/dahlcurrent.pdf .  Cite
as  arXiv:1110.4678 . 

\item{[EW]}J.~Eisert and M.~Wilkens, ``Quantum Games'',
{\it J.
of Modern Optics\/} 47 (2000), 2543-2556

\item{[EWL]} J. Eisert, M. Wilkens and
M. Lewenstein, ``Quantum Games and
Quantum Strategies'', {\it Phys. Rev.
Letters\/} 83, 3077 (1999).

\item{[K]}H. Kuhn, ``Extensive Games and the Problem of
Information'', in {\it Contributions to the Theory of
Games
II}, Annals of Math Studies 28 (1953).

\item{[L]}S.~Landsburg, ``Nash Equilibria in Quantum
Games'', {\it Proc. Am. Math. Soc.\/}, December 2011. 

\item{[M]}D.~Meyer, ``Quantum Strategies'', {\it Phys.
Rev.
Lett.\/} 82 (1999), 1052-1055

\bye